\documentclass{article}
\usepackage{amssymb}
\usepackage{amsfonts}
\usepackage{amsmath}
\usepackage{geometry}
\usepackage{graphicx}

\begin{document}

\title{Incremental Similarity and Turbulence} 
\author{Ole E. Barndorff-Nielsen, Emil Hedevang and J\"urgen Schmiegel\\
Department Of Mathematics, Aarhus University, Denmark}
\date{}
\maketitle

\begin{abstract}
This paper discusses the mathematical representation of an empirically observed phenomenon, referred to as Incremental Similarity. We discuss this feature from the viewpoint of stochastic processes and present a variety of non-trivial examples, including those that are of relevance for turbulence modelling.
\end{abstract}

\section{Introduction}

This paper discusses an empirically observed feature that we refer to as 
\emph{Incremental Similarity}. Its most remarkable manifestation is in
turbulence but it is also found, in a less prononunced form, in finance.

In high frequency recordings of velocities of homogeneous, isotropic and
stationary turbulence the non-Gaussian and skewed distributions of velocity
increments from different experiments and different lags are essentially
identical provided they have the same variance. Note however that this
identification is effectuated lag by lag and not by a simple common
transformation.

The observation of incremental similarity for turbulent velocity time series
adds a new type of universality to the statistical stylised features of
turbulent flows. In the time domain, the term universality traditionally
refers to universal scaling properties of structure functions defined as
moments of velocity increments. These scaling laws are only realised in the
limit of very large Reynolds numbers, i.e. strongly turbulent flows. Even
for large Reynolds numbers only approximate scaling is observed and only for
a certain range of time scales, the so-called inertial range, which covers
merely part of the dynamically active scales.

The new stylised feature of incremental similarity points towards a
completely different type of universality. Incremental similarity
characterises the distributions of velocity increments at all dynamically
active scales. Furthermore, our empirical analysis shows that incremental
similarity is not restricted to any high Reynolds number limit. And finally,
incremental similarity provides a relatively simple mapping that directly
connects measurements from different experiments, different in Reynolds
numbers and different in boundary conditions.

The phenomenon of incremental similarity was first noted in the paper \cite%
{BNBlSch04} which revealed this trait through a detailed analysis of the
recordings of the main component of the velocity vector from three different
types of experiments, with Reynold's numbers 80, 190, 17000. The
collapsibility of the incremental distributions was of extraordinary degree.
This type of analysis has since been extended to a much larger class of
experiments, as reported in \cite{BNHedSch15}, fully confirming the original
observation. However, up till now a mathematical representation of the
phenomenon has not been given. Such a representation is proposed in this
paper.

As stated here the property of incremental (or distributional) similarity is
of nonparametric character but, as demonstrated in \cite{BNBlSch04}, the
laws of the velocity differences can be fitted by the normal inverse
Gaussian distribution, denoted $NIG$, to a high degree of precision. This
distribution was introduced in \cite{BN79} and has since found a multitude
of applications in a variety of fields. Recently the fact that $NIG$
describes the velocity increments to high precision has found a
theoretical counterpart in studies of Birnir concerning a stochastic version
of the Navier-Stokes equations, see \cite{Bir13a,Bir13b,Bir14}.

The empirical evidence for incremental similarity is summarised in Section 2. Section 3 presents a mathematical definition of incremental similarity - or IS - and an extended concept - extIS - and a variety of examples of stochastic processes and fields that
meet the definition of incremental similarity exactly are presented there. For the relation to the statistical theory of turbulence it is crucial to have examples where the processes considered are stationary, and herein lies the main difficulty. Section 4 briefly adresses the question of modelling the timewise behaviour
of the main component of the velocity vector in homogeneous turbulence by
stochastic processes embodying the main universal features of this type of
dynamics, including that of incremental similarity. 

\section{Incremental similarity: empirical evidence}

The empirical verification of incremental similarity, presented in \cite{BNBlSch04}\ and \cite{BNHedSch15} and briefly outlined here, is based on
the analysis of 17 experimental data sets (a wake-flow experiments, a
free-jet experiment, 13 helium jet experiments, one wind tunnel experiment,
one data set from the atmospheric boundary layer) and on one data set from a
direct numerical simulation (DNS) of the Navier-Stokes equation. The
Reynolds numbers covered by these data range from 80 up to 20000. The
empirical data consist of stationary time series of recordings of the main
component of the turbulent velocity vector measured at a fixed position in
space. The DNS data are spatially homogeneous at a fixed point in time. All
data are standardised for the velocities to have a unit variance.

Let $v_{t}^{(i)}$ denote the velocity signal at time $t$ belonging to the
data set $(i)$. We denote by 
\begin{equation*}
u^{(i)}_{s}=v^{(i)}_{t+s}-v^{(i)}_{t}
\end{equation*}
the velocity increment at time scale $s$. Here we skip reference to $t$
since we are only dealing with stationary time series. (For the DNS data $t$
denotes the spatial position.)

\begin{figure}
  \centering
  \includegraphics{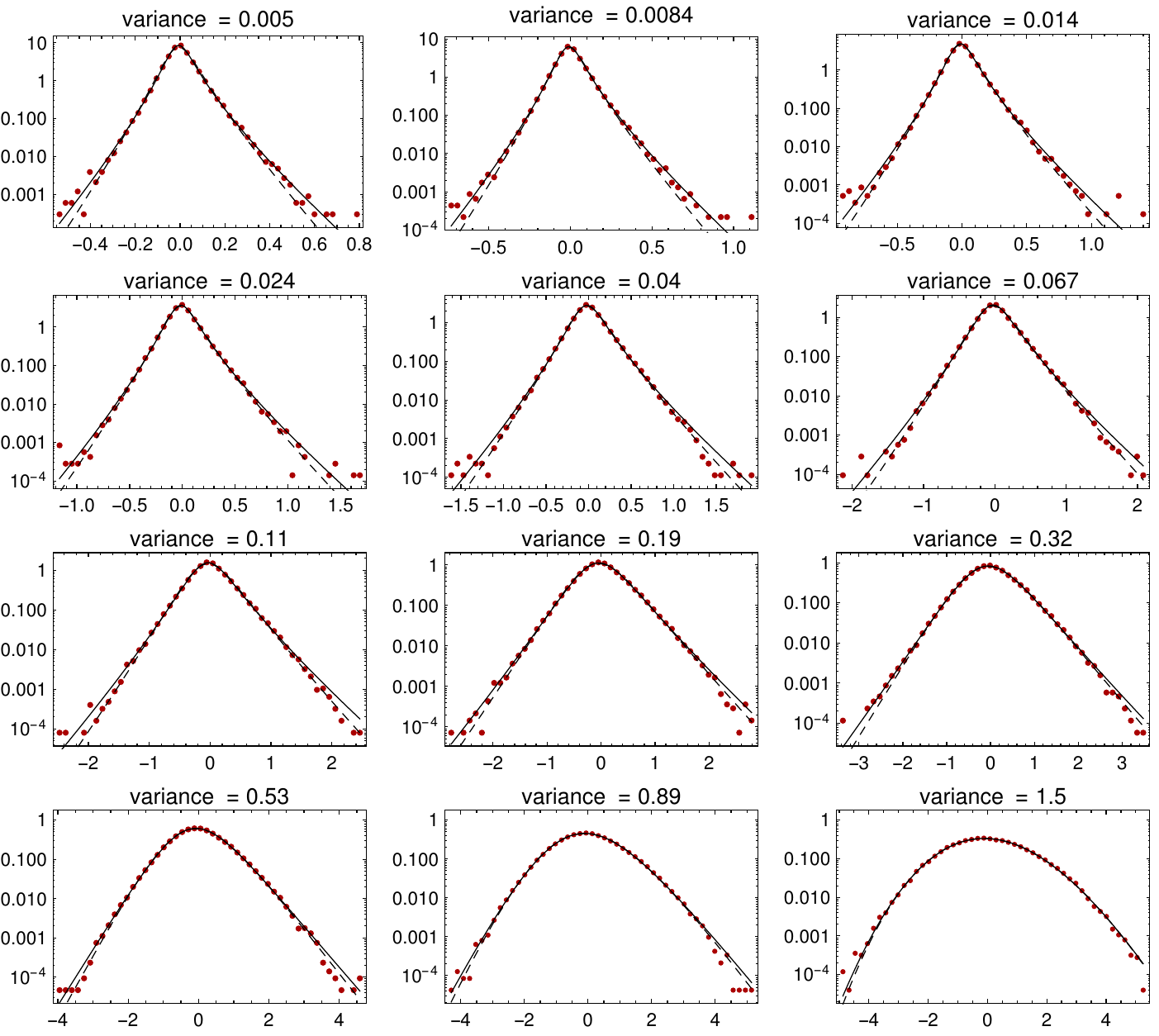}
  \caption{Logarithmic representation of probability density
    functions. Red dots: empirical pdf of the wind tunnel dataset.
    Solid black: pdf of the corresponding estimated normal inverse
    Gaussian distribution. Dashed black: pdf of the normal inverse
    Gaussian distribution fitted from all sixteen data sets pooled
    into a single one.}
  \label{fig:nig-fit}
\end{figure}

A key observation related to the densities of velocity increments is
depicted in Figure \ref{fig:nig-fit}. Each graph corresponds to the densitiy
of velocity increments at a certain time scale $s$ with $s$ increasing from
top left to bottom right. These densities evolve from heavy tails at small
time scales towards a more Gaussian shape at the large time scales. This
evolution across scales is well known in the literature and sometimes called 
aggregational Gaussianity. Figure \ref{fig:nig-fit} also shows, as
solid lines, the approximation of these densities within the class of normal
inverse Gaussian distributions. The normal inverse Gaussian distributions
fit the empirical densities equally well for all time lags and all
amplitudes.

Figure \ref{fig:nig-fit} provides one example for the evolution across
scales, similar results are observed for the other data sets we analysed.
The larger the Reynolds number the more the heavyness of the tails. And,
different data sets show different distributions at the same time scale. But
the key question that is of interest here is: Do different experiments show
the same distribution of increments just at different time scales? In other
words: Is the evolution across scales exemplified in Figure \ref{fig:nig-fit}
universal in the sense that the distributions of $u_{s}^{(i)}$ and $%
u_{s^{\prime }}^{(j)}$ are the same if $s^{\prime }$ is properly chosen
given $s$? Obviously, an affirmative answer requires that the variances at
these time scales are the same 
\begin{equation*}
\mathrm{Var}\left( u_{s}^{(i)}\right) =\mathrm{Var}\left( u_{s^{\prime
}}^{(j)}\right) .
\end{equation*}

\begin{figure}[ht]
  \centering
  \includegraphics{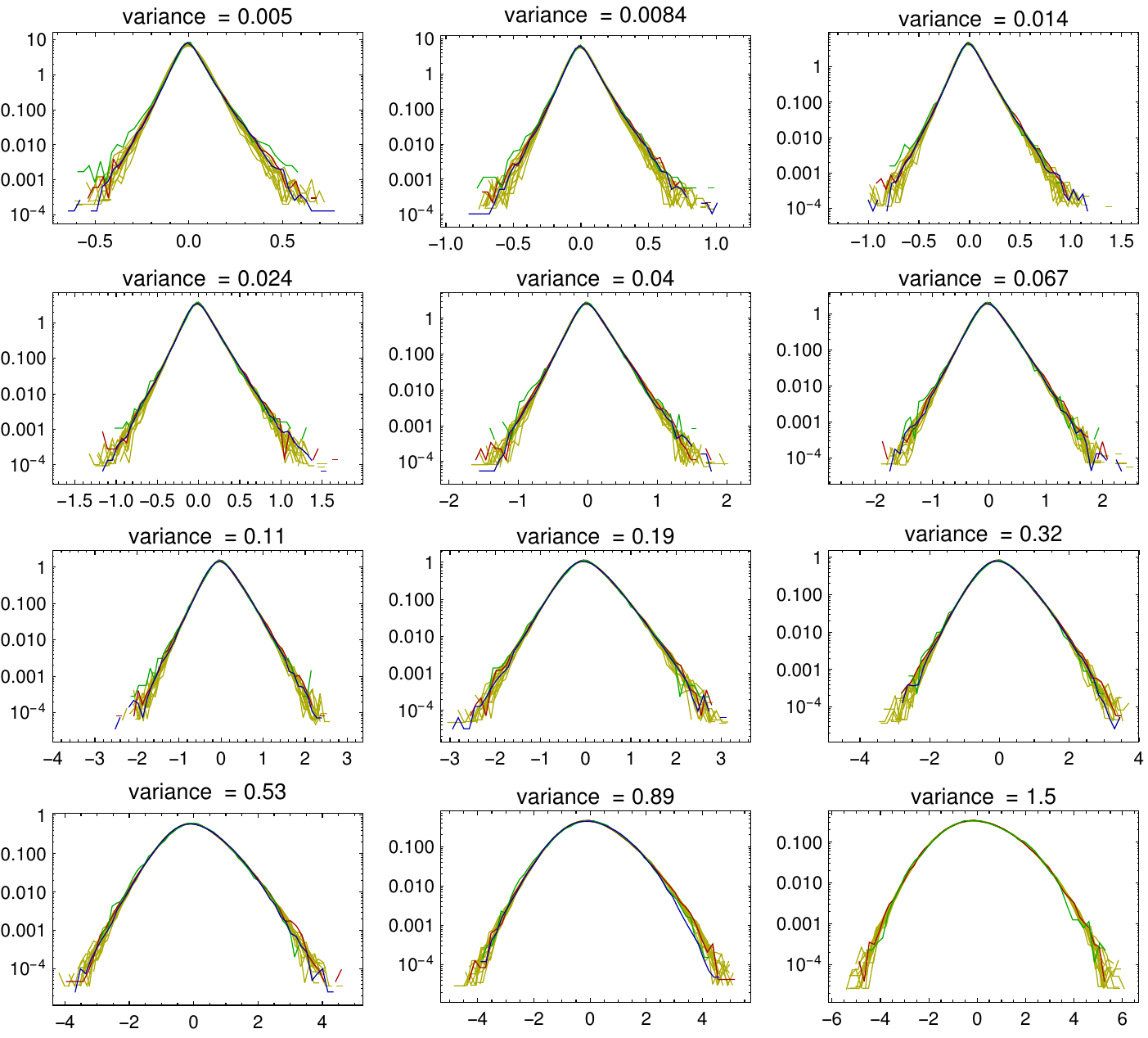}
  \caption{Logarithm representation of the empirical probability
    density function of the sixteen datasets, after pairing of variances. The thirteen Helium jet data
    sets are displayed using the same color (yellow).The wind tunnel dataset is displayed in red and the data set from the atmospheric boundary layer is displayed in green. The DNS dataset (blue) is
    absent in the case of the largest variance.}
  \label{fig:inc-same-var}
\end{figure}

Figure \ref{fig:inc-same-var} shows the corresponding densities of velocity
increments for $12$ fixed values of the variance. Each plot corresponds to a
different value, increasing from top left to bottom right. The densities
within each plot correspond to different experiments and different time
lags, but the variances are the same. We clearly observe the collapse of the
densities during the whole evolution across scales. In this sense the
evolution across scales is universal.

The shape of the distributions in Figure \ref{fig:inc-same-var} is, to a
good approximation, a universal function of the variance. For normal
distributions, this is trivial. But here we clearly have distributions that
are not normal.

\section{Incremental similarity: mathematical considerations}

As theoretical counterpart to the incremental similarity features discussed
in Sections 1 and 2 we introduce the following definition\newline

\textbf{Definition}\quad \emph{\ Incremental similarity}\quad Let $X$\ and $Z
$\ be two stochastic processes on $\mathbb{R}$. Then $X$\ and $Z$\ are said
to be incrementally similar, or IS, provided that for any $t\in \mathbb{R}$
and any $u>0$\ there exists a $t^{\prime }\in \mathbb{R}$\ and a positive
number $u^{\prime }$\ such that the law of $Z\left( t^{\prime }+u^{\prime
}\right) -Z\left( t^{\prime }\right) $ is the same as that of $X\left(
t+u\right) -X\left( t\right) $.

More generally, if $\mathcal{X}$\ and $\mathcal{Z}$ are two classes of
stochastic process on $\mathbb{R}$ then $\mathcal{X}$\ and $\mathcal{Z}$\
are said to have the IS\ property if all pairs $X$\ and $Z$ such that $X\in 
\mathcal{X}$\ and $Z\in \mathcal{Z}$ are of IS type. For brevity we will
then say that $\left( \mathcal{X},\mathcal{Z}\right) $ is IS and that $%
\mathcal{X}$\ is IS if that is the case of $\left( \mathcal{X},\mathcal{X}%
\right) $.\quad $\square $\newline

There are many trivial examples of increment similarity. For instance, any
continuous Gaussian process is incrementally similar to Brownian motion.
Another example is where $X$ and $Z$\ \ are stationary processes on $\mathbb{%
R}_{+}$\ with $Z$\ equal in law to a proportional timechange of $X$, that is 
$Z_{t}=X_{ct}$\ for some $c>0$. Various non-obvious examples may be based on
the following concept of extended incremental similarity or extIS.\newline

\textbf{Definition}\quad \emph{\ Extended incremental similarity}\quad Let $%
R $\ be a class of positive, continuous and decreasing functions $r$\ on $%
[0,\infty )$. Further, let $\mathcal{X}=\left\{ X^{\left[ r\right] }:r\in
R\right\} $\ be a parametrised family of stationary processses $X^{\left[ r%
\right] }$ on $\mathbb{R}$ with the property that for any pair of time
points $\left( t,t+u\right) $\ the joint law of $X_{t}^{\left[ r\right] }$\
and $X_{t+u}^{\left[ r\right] }$\ is fully determined by $r\left(u\right) $%
; this in particular implies that the same holds for the law of the
increment $X_{t+u}^{\left[ r\right] }-X_{t}^{\left[ r\right] }$. We
furthermore assume that the joint law of $X_{t}^{\left[ r\right] }$\ and $%
X_{t+u}^{\left[ r\right] }$ is the same as the joint law of $X_{t}^{\left[ 
\tilde{r}\right] }$\ and $X_{t+\tilde{u}}^{\left[ \tilde{r}\right] }$\
provided $r\left( u\right) =\tilde{r}\left( \tilde{u}\right) $ whatever $r$\
and $\tilde{r}$\ in $R$. We then say that, relative to $R$, $\mathcal{X}$\
has the property of extended incremental similarity or that $\mathcal{X}$\
is extIS. \quad $\square $\newline

If $\mathcal{X}$\ is extIS then it is in particular IS. As a another direct
consequence of the definition of extIS we have\newline

\textbf{Proposition}\quad Let $\mathcal{X}_{1},\mathcal{X}_{2},...,\mathcal{X%
}_{n}$ denote independent extIS families relative to the same index class $R$%
. If $F$\ is a real (measurable) function on $\mathbb{R}^{n}$\ and if $%
\mathcal{Y}=\left\{ Y^{\left[ r\right] }:r\in R\right\} $\ is the class of
processes given by%
\begin{equation*}
Y_{t}^{\left[ r\right] }=F\left( X_{1t}^{\left[ r\right] },...,X_{nt}^{\left[
r\right] }\right) 
\end{equation*}%
where $X_{j}^{\left[ r\right] }\in \mathcal{X}_{j}$, $j=1,...,n$,\ then the
family $\mathcal{Y}$\ is extIS, and hence IS. In fact, the same conclusion
holds if $F$\ is random, provided it is independent of $\left( \mathcal{X}%
_{1},\mathcal{X}_{2},...,\mathcal{X}_{n}\right) $.\quad $\square $\newline

Applications of this result will be discussed in Section 4.

We proceed to present some classes of stationary processes on $\mathbb{R}$
having the extIS property

Let $\mathcal{U}$ be the class\ of stationary Gaussian processes on $\mathbb{%
R}$ of mean $0$\ and variance $1$ such that for any member $X$\ of $\mathcal{%
U}$\ the autocorrelation function $r\ $of $X$\ is positive, continuous and
strictly decreasing to $0$. Then, as is easily seen, $\mathcal{U}$\ is extIS.

The concept of trawl processes, introduced in \cite{BN11}, offers a range of
extIS classes. As discussed in \cite{BNBV15}, the simplest type of trawl
processes $X$\ on $\mathbb{R}$\ are of the form%
\begin{equation}
X_{t}=L\left( A_{t}\right)   \label{Xtrawl}
\end{equation}%
where $L$\ is a homogeneous L\'{e}vy basis on $\mathbb{R}^{2}$ and $%
A_{t}=A+\left( t,0\right) $\ for a Borel set $A$ in $\mathbb{R}^{2}$ with
positive Lebesgue measure, points in $\mathbb{R}^{2}$\ being denoted by $%
\left( t,x\right) $. In this case, since for any $\phi ,\psi \in \mathbb{R}$
we have%
\begin{eqnarray*}
\phi X_{t}+i\psi X_{t+u} &=&\phi L\left( A_{t}\backslash A_{t+u}\right)
+\left( \phi +\psi \right) L\left( A_{t}\cap A_{t+u}\right) +\psi L\left(
A_{t+u}\backslash A_{t}\right)  \\
&=&\phi L\left( A\backslash A_{u}\right) +\left( \phi +\psi \right) L\left(
A\cap A_{u}\right) +\psi L\left( A_{u}\backslash A\right) 
\end{eqnarray*}%
the cumulant function of $\left( X_{t},X_{t+u}\right) $ is given by%
\begin{eqnarray*}
\mathrm{C}\{\phi ,\psi \ddag \left( X_{t},X_{t+u}\right) \} &=&\left\vert
A\backslash A_{u}\right\vert \mathrm{C}\{\phi \ddag L^{\prime }\}+\left\vert
A\cap A_{u}\right\vert \mathrm{C}\{\phi +\psi \ddag ^{\prime }\}+\left\vert
A_{u}\backslash A\right\vert \mathrm{C}\{\phi \ddag L^{\prime }\} \\
&=&\left\vert A\right\vert \left[ \mathrm{C}\{\phi \ddag L^{\prime }\}+%
\mathrm{C}\{\psi \ddag L^{\prime }\}+r\left( u\right) \left( \mathrm{C}%
\{\phi \ddag L^{\prime }\}+\mathrm{C}\{\psi \ddag L^{\prime }\}-\mathrm{C}%
\{\phi +\psi \ddag ^{\prime }\}\right) \right] 
\end{eqnarray*}%
where $L^{\prime }$\ denotes the L\'{e}vy seed of $L$, $\left\vert
{\cdot}\right\vert $ indicates Lebesgue measure and%
\begin{equation}
r\left( u\right) =\frac{\left\vert A\cap A_{u}\right\vert }{\left\vert
A\right\vert }.  \label{trawl r}
\end{equation}

Now consider the class $\mathcal{A}$\ of Borel sets $A$ such that (\ref%
{trawl r})\ is positive for all\ real $u$ with $r$\ continuous and strictly
decreasing on $\mathbb{R}_{+}$ and tending to $0$\ as $u\longrightarrow
\infty $. Let $R$ be the corresponding class of functions $r$. Furthermore,
for any $c>0$\ let $\mathcal{A}_{c}$\ be the subclass of $\mathcal{A}$\
given by $\mathcal{A}_{c}=\left\{ A\in \mathcal{A}:\left\vert A\right\vert
=c\right\} $. If $\mathcal{X}$\ is the class of trawl processes
corresponding to a given seed $L^{\prime }$\ and a given $\mathcal{A}_{c}$
then $\mathcal{X}$\ is extIS and hence IS.

We note that%
\begin{equation*}
\mathrm{C}\{\phi \ddag X_{t+u}-X_{t}\}=\left\vert A\right\vert \left(
1-r\left( u\right) \right) \left( \mathrm{C}\{\phi \ddag L^{\prime }\}+%
\mathrm{C}\{-\phi \ddag L^{\prime }\}\right) \mathrm{.}
\end{equation*}

In case $L^{\prime }$\ is square integrable then $r$\ is the autocorrelation
function of the process $Y$. In general we will refer to $r$\ as the \emph{%
autodependence} function of $Y$. By suitable choice of $A$\ the
autodependence function can be selected to show short, middle or long term
dependence.

With reference to Section 2 we note that if the law of $L^{\prime }$\ is
normal inverse Gaussian and symmetric then $Y_{t}$\ and all increments of $Y$%
\ are also normal inverse Gaussian distributed.

Next, suppose that $X_{t}\ $\ is a stationary process of the form%
\begin{equation}
X_{t}=\int_{-\infty }^{t}g\left( t-s\right) L\left( \mathrm{d}s\right)
\label{LSS}
\end{equation}%
where $L$\ is a symmetric $\alpha $-stable L\'{e}vy process on $\mathbb{R}$\
and the kernel function $g$ satisfies $g\left( s\right) =0$\ for $s<0$\ and%
\begin{equation}
I_{\alpha }\left( g\right) =\int_{0}^{\infty }g\left( s\right) ^{\alpha }%
\mathrm{d}s<\infty .  \label{Ia}
\end{equation}%
for some $\alpha \in \left( 0,2\right) $. Then the integral (\ref{LSS})\
exists and, since the L\'{e}vy seed $L^{\prime }$\ of $L$ has cumulant
function%
\begin{equation}
\mathrm{C}\{\phi \ddag X\}=-\gamma \left\vert \phi \right\vert ^{\alpha },
\end{equation}%
we find, for $u>0$, that%
\begin{equation}
\mathrm{C}\{\phi ,\psi \ddag X_{t},X_{t+u}\}=-\gamma \left\vert \phi
\right\vert ^{\alpha }\int_{0}^{u}g\left( s\right) ^{\alpha }\mathrm{d}%
s-\gamma \int_{-\infty }^{0}\left\vert \phi g\left( u-s\right) -\psi g\left(
-s\right) \right\vert ^{\alpha }\mathrm{d}s.  \label{cumaXX}
\end{equation}%
Now, let $\mathcal{G}$\ be the class of kernel functions $g$\ such that $g$\
is continuous and strictly decreasing to $0$ and let $\mathcal{X}$\ be the
corresponding class of stochastic processes (\ref{LSS}). Suppose that $g$\
and $h$ are both members of $\mathcal{G}$\ and consider the analogue of (\ref%
{cumaXX}), i.e.%
\begin{equation}
\mathrm{C}\{\phi ,\psi \ddag Z_{t},Z_{t+u}\}=-\gamma \left\vert \phi
\right\vert ^{\alpha }\int_{0}^{u}h\left( s\right) ^{\alpha }\mathrm{d}%
s-\gamma \int_{-\infty }^{0}\left\vert \phi h\left( u-s\right) -\psi h\left(
-s\right) \right\vert ^{\alpha }\mathrm{d}s
\end{equation}%
where $Z$\ denotes the element of $\mathcal{X}$\ corresponding to $h$. Only
in quite exceptional cases will it be possible for every $u>0$\ to find a $%
v>0$\ such that $\mathrm{C}\{\phi ,\psi \ddag Z_{t},Z_{t+v}\}=\mathrm{C}%
\{\phi ,\psi \ddag X_{t},X_{t+u}\}$. In other words, in the present setting
interesting examples of extIS do not exist.

On the other hand,%
\begin{equation}
\mathrm{C}\{\phi \ddag X_{t+u}-X_{t}\}=-\gamma \left\vert \phi \right\vert
^{\alpha }\hat{g}\left( u;\alpha \right)   \label{cumaDX}
\end{equation}%
where%
\begin{equation}
\hat{g}\left( u;\alpha \right) =\int_{-\infty }^{u}\left\vert g\left(
u-s\right) -g\left( -s\right) \right\vert ^{\alpha }\mathrm{d}s  \label{ghat}
\end{equation}%
and therefore there are subclasses of $\mathcal{X}$ that are IS.
Specifically, for a fixed $\alpha \in \left( 0,2\right) $, consider the
subclass\ $\mathcal{G}_{\alpha }$\ of $\mathcal{G}$ of kernels $g$\ such
that $I_{\alpha }\left( g\right) $ does not depend on $g$. Then the
processes $X$\ in\ $\mathcal{X}$\ have the same one-dimensional marginal
distribution, and $\mathcal{X}$\ is IS. In fact, for any $g,h\in \mathcal{G}%
_{\alpha }$\ and any $u>0$ there exists a $v>0$\ such that the condition $%
\left\vert \hat{h}\left( v;\alpha \right) \right\vert =\left\vert \hat{g}%
\left( u;\alpha \right) \right\vert $ is met.

From the viewpoint of applications the question now is whether in principle
it is possible, given a class of processes that are known to be of IS type
(or suspected to be so), to determine the transformation that effectuates
the collapsibility of the laws of increments from the different members of
the class.

Empirically, given that high frequency and extensive datasets are available
from each of the processes in question and that stationarity of the series
is a realistic assumption, the most immediate way is to first standardise
the series to have the same marginal variance and then estimate the
variances of the increments for a suitable range of lags, as was done for
the turbulence data discussed in Section 2. However, in case the data are
suspected to come from a fractional regime where second order moments may
not in principle exist, such as the $\alpha $-stable setting considered
above, one may resort to other methods of lag transformation, for instance
resorting to estimation of fractional moments.

\section{Models of $\mathcal{BSS}$/$\mathcal{LSS}$ type}

In \cite{BNSch09}\ the concept of Brownian semistationary processes - or $%
\mathcal{BSS}$\ - processes was introduced. These are stationary processes
of the form%
\begin{equation}
Y_{t}=\mu +\int_{-\infty }^{t}g\left( t-s\right) \sigma _{s}\mathrm{d}%
B_{s}+\int_{-\infty }^{t}q\left( t-s\right) a_{s}\mathrm{d}s  \label{BSS}
\end{equation}%
where $B$\ is Brownian motion, $\sigma $ and $a$\ are stationary processes
and the kernels $g$ and $q$ are deterministic functions. The particular
setting%
\begin{equation}
Y_{t}=\mu +\int_{-\infty }^{t}g\left( t-s\right) \sigma _{s}\mathrm{d}%
B_{s}+\int_{-\infty }^{t}q\left( t-s\right) \sigma _{s}^{2}\mathrm{d}s
\label{BSSpar}
\end{equation}%
is of particular interest. (It may be seen as a stationary analogue of the
BNS\ model studied in financial econometrics.) The $\mathcal{LSS}$ class is
obtained by substitution of $B$\ in (\ref{BSS}) by a general L\'{e}vy
process $L$.

The primary aim of the definitions of $\mathcal{BSS}$\ and $\mathcal{LSS}$
was to model the timewise behaviour of the main velocity component in a
homogeneous turbulent flow, but the same kind of processes have found many
applications elsewhere, see for instance \cite{BNBV15}, \cite{Pod14} and
references given there.

$\mathcal{BSS}$ processes of type (\ref{BSSpar})\ have been demonstrated to
be capable of modelling classical stylised features of turbulence very
accurately, cf. \cite{MarSch15}. And a recent study \cite{HedSch14} of the
Helium data discussed in Section 2, based on exponential $\mathcal{LSS}$
processes, has revealed a new type of universality for the energy dissipation.

However, the $\mathcal{BSS}$\ structure does not have the property of
extended increment similarity and it is therefore natural to ask whether
there is a modification of that structure having dynamic behaviour closely
similar to $\mathcal{BSS}$\ but also exhibiting extIS.

To this end we now introduce the following variant of $\mathcal{BSS}$. Let $%
X $\ be a stationary Gaussian process such that the\ one-dimensional
marginals of $X$\ have mean $0$\ and variance $1$. Let $\sigma $ be a
stationary, cadlag and square integrable process with autocorrelation
function $\rho $, and let 
\begin{equation}
Y_{t}=\mu +\sigma _{t}X_{t}+\beta \sigma _{t}^{2}  \label{BSS'}
\end{equation}%
for some constant $\beta \in \mathbb{R}$. We shall refer to this type as a $%
\mathcal{BSS}^{\prime }$\ process. Colloquially speaking, the alternative
view amounts to moving the volatility process $\sigma $ in (\ref{BSSpar})
outside the integration signs.

Now, let $\mathcal{Y}$\ be the class of processes (\ref{BSS'}) obtained by
letting the autocorrelation function $r$ of $X$\ vary over the set $R$\ for
which $r$\ is positive, continous and decreasing on $\mathbb{R}_{+}$ and
taking $\sigma ^{2}$\ to vary over a family $\Sigma $\ which in itself is
extIS with the same index set $R$ and such that $\sigma $\ is paired with $X$%
\ so as to have the same index $r$. Then, in view of the Proposition
presented in Section 3, the class $\mathcal{Y}$\ is extIS.

For instance, $\Sigma $\ could be chosen to be an extIS class of processes $%
\sigma ^{2}$\ for which the $\log \sigma ^{2}$\ are trawl process $L\left(
A+\left( t,0\right) \right) $, as discussed in Section 3 \ In particular, $L$
may be taken to be an $NIG$\ basis, as in \cite{HedSch14}.

As a further example, suppose that $Z$\ is an independent copy of $X$\ and
let $\sigma _{t}^{2}=\left\vert Z_{t}\right\vert ^{1/2}$. In this case the
autovariance function of the corresponding process $Y$\ of (\ref{BSS'})\ is
given by 
\begin{equation*}
\mathrm{E}\left\{ \left( Y_{t}-Y_{0}\right) ^{2}\right\} =2\left[ \mathrm{E}%
\left\{ \sigma _{0}\right\} ^{2}\left( 1-r\left( t\right) \right) +\mathrm{V}%
\left\{ \sigma _{0}\right\} r\left( t\right) \left( 1-\rho \left( t\right)
\right) +2\beta ^{2}\mathrm{V}\left\{ \sigma _{0}^{2}\right\} \left(
1-\varrho \left( t\right) \right) \right] 
\end{equation*}%
where $\varrho $\ is the autocorrelation function of $\sigma ^{2}$, with%
\begin{equation*}
\mathrm{E}\left\{ \left( Y_{t}-Y_{0}\right) ^{2}\right\} =2[\left( 3+96\beta
^{2}\right) -r\left( t\right) -2r\left( t\right) ^{3}-72\beta ^{2}r\left(
t\right) ^{2}-24\beta ^{2}r\left( t\right) ^{4}
\end{equation*}%
\ which is a monotonely decreasing function of $r\left( t\right) $ enabling
direct lag identification.

\section*{Acknowledgements}
The authors are much indebted to K.R. Sreenivasan, J. Peinke, B. Chabaud and C. Meneveau for allowing to use the data sets.

\newpage

\end{document}